\documentstyle[a4,11pt]{article}
\begin{document}
\begin{center}
ON THE SOLITON GEOMETRY IN MULTIDIMENSIONS
\end{center}
\vskip 2cm
\begin{center}
(Talk prepared to the conference "Differential Geometry
and Quantum Physics", Berlin, 06-10.03.2000)
\end{center}
\vskip 2cm
\begin{center}
R.N.Syzdykova\footnote{E-mail: cnlprns@satsun.sci.kz},
Kur. R. Myrzakul\footnote{Permanent address:
Department of Mathematics and Mechanics, Kazakh State University},
G.N.Nugmanova\footnote{E-mail: cnlpgnn@satsun.sci.kz}
and R.Myrzakulov\footnote{E-mail: cnlpmyra@satsun.sci.kz}
\end{center}
\vskip 1cm

\begin{center}
 Institute of Physics and Technology, 480082 Alma-Ata-82, Kazakhstan
\end{center}

\begin{abstract}

The connection between multidimensional soliton equations and differential
geometry of surfaces/curves is discussed. A classes integrable surfaces
(and/or curves) corresponding to the well known (2+1) - dimensional soliton
equations are constructed. Examples of such equations include: the Davey -
Stewartson, Zakharovs, Kadomtsev - Petviashvili, (2+1)-mKdV, (2+1)-KdV
and Knizhnik - Novikov - Veselov equations.    Using the presented geometrical formalism
the equivalence between these soliton equations and their spin  counterparts
are established. The integrable curve/surface corresponding to the self-dual
Yang-Mills equation is constructed.

\end{abstract}

\tableofcontents
\section{Introduction}
Recently many works have been devoted to the study on motion of curves
ad surfaces. One of the interesting motivations is the relation to soliton
equations [1-10]. In particular, the
geometrical formulation may set a suitable basis for the description
of higher dimensional soliton equations [9-14].
In this note we would like to discuss the simplest and most
transparent differential-geometric appearance of the some soliton
equations in multidimensions (see, also , [16-27]).

Three equations of the differential geometry are may be starting points
of the soliton geometry: the Serret-Frenet equation (SFE),  the
Gauss-Weingarten equation (GWE) and the
Gauss-Mainardi-Codazzi equation (GMCE). We begin with the presentation
of these equations.

\subsection{The Serret-Frenet equation}
The SFE has the form
$$
\left ( \begin{array}{ccc}
{\bf e}_{1} \\
{\bf e}_{2} \\
{\bf e}_{3}
\end{array} \right)_{x}=
\left ( \begin{array}{ccc}
0             & k     &  0 \\
-\beta k      & 0     & \tau  \\
0             & -\tau & 0
\end{array} \right)
\eqno(1)
$$
where $k, \tau $ are the curvature and torsion, $x$ is the arc-length
parameter, ${\bf e}_{1}, {\bf e}_{2} $ and ${\bf e}_{3}$ are the tangent,
normal and binormal vectors, respectively with
$$
{\bf e}_{1}^{2}=\beta = \pm 1, {\bf e}_{2}^{2}={\bf e}^{2}_{3}=1 \eqno(2)
$$
where $\beta=+1$ corresponds the the focusing case and $\beta=-1$ corresponds
to the defocusing case.

\subsection{The Gauss-Weingarten equation}
Let us consider a surface in $R^{3}$ equipped (locally) with coordinates
$x, y$ and defined by the position vector ${\bf r}={\bf r}(x,y) \in R^{3}$.
Then the GWE reads as
$$
{\bf r}_{xx} = \Gamma^{1}_{11} {\bf r}_{x} + \Gamma^{2}_{11} {\bf r}_{y}
+ L{\bf n}   \eqno(3a)
$$
$$
{\bf r}_{xy} = \Gamma^{1}_{12} {\bf r}_{x} + \Gamma^{2}_{12} {\bf r}_{y}
+ M{\bf n}   \eqno(3b)
$$
$$
{\bf r}_{yy} = \Gamma^{1}_{22} {\bf r}_{x} + \Gamma^{2}_{22} {\bf r}_{y}
+ N{\bf n}   \eqno(3c)
$$
$$
{\bf n}_{x} = p_{11} {\bf r}_{x} + p_{12} {\bf r}_{y}   \eqno(3d)
$$
$$
{\bf n}_{y} = p_{21} {\bf r}_{x} + p_{22} {\bf r}_{y}   \eqno(3e)
$$
where $\Gamma^{k}_{ij}$ are the Christoffel symbols. This equation can be
rewritten in the form
$$
Z_{x} = AZ, \quad Z_{y} = BZ  \eqno(4)
$$
where
$Z=({\bf r}_{x}, {\bf r}_{y}, {\bf n})^{t}$ and
$$
A =
\left ( \begin{array}{ccc}
\Gamma^{1}_{11} & \Gamma^{2}_{11} & L \\
\Gamma^{1}_{12} & \Gamma^{2}_{12} & M \\
p_{11}           & p_{12}           & 0
\end{array} \right) ,  \quad
B =
\left ( \begin{array}{ccc}
\Gamma^{1}_{12} & \Gamma^{2}_{12} & M \\
\Gamma^{1}_{22} & \Gamma^{2}_{22} & N \\
p_{21}          & p_{22}          & 0
\end{array} \right). \eqno(5)
$$

On the other hand, in terms of the orthogonal frame
$$
{\bf e}_{1} = \frac{{\bf r}_{x}}{\sqrt E}, \quad
{\bf e}_{2} = {\bf n},
\quad {\bf  e}_{3} = {\bf e}_{1} \wedge {\bf  e}_{2}    \eqno(6)
$$
the GWE takes the form
$$
\left ( \begin{array}{ccc}
{\bf e}_{1} \\
{\bf e}_{2} \\
{\bf e}_{3}
\end{array} \right)_{x}=A
\left ( \begin{array}{c}
{\bf e}_{1} \\
{\bf e}_{2} \\
{\bf e}_{3}
\end{array} \right), \quad
\left ( \begin{array}{ccc}
{\bf e}_{1} \\
{\bf e}_{2} \\
{\bf e}_{3}
\end{array} \right)_{y}=B
\left ( \begin{array}{c}
{\bf e}_{1} \\
{\bf e}_{2} \\
{\bf e}_{3}
\end{array} \right).
\eqno(7)
$$
Here
$$
A=\frac{1}{\sqrt{E}}
\left ( \begin{array}{ccc}
0             & L   & -\frac{g}{\sqrt{E}}\Gamma^{2}_{11} \\
-\beta L      & 0     & -g p_{12}  \\
\frac{\beta g}{\sqrt{E}}\Gamma^{2}_{11}& g p_{12} & 0
\end{array} \right), \quad
B=\frac{1}{\sqrt{E}}
\left ( \begin{array}{ccc}
0             & M   & -\frac{g}{\sqrt{E}}\Gamma^{2}_{12} \\
-\beta M      & 0     & -g p_{22}  \\
\frac{\beta g}{\sqrt{E}}\Gamma^{2}_{12}& g p_{22} & 0
\end{array} \right) \eqno(8)
$$
where $g=EG-F^{2}.$

\subsection{The Gauss-Mainardi-Codazzi equation}

Integrability conditions for the GWE can be derived from
$$
{\bf r}_{xxy}= {\bf r}_{xyx}, \quad {\bf r}_{xyy}= {\bf r}_{yyx} \eqno(9a)
$$
or
$$
{\bf e}_{jxy}= {\bf e}_{jyx} \eqno(9b)
$$
and are equivalent to the following GMCE
$$
A_{y}-B_{x}+[A,B]=0 \eqno(10)
$$
where $A, B$ are given by (5) or (8). In general, this equation is apparently nonintegrable.
Only in some particular cases
it is integrable (see, e.g. [33]). Also it admits some multidimensional extensions,
e.g. in 2+1  dimensions. It is remarkable that some of these
(2+1)-dimensional extensions are integrable, e.g. the M-LXII equation (19).

\section{Soliton geometry in $d=2$ dimensions}

To make exposition self-contained, we expose here
some neccessary to us well known basis facts from the 2-dimensional
soliton geometry. We must construct two-dimensional generalization
of the SFE (1).  The simplest and well known two-dimensional extension of the SFE has the form
$$
\left ( \begin{array}{ccc}
{\bf e}_{1} \\
{\bf e}_{2} \\
{\bf e}_{3}
\end{array} \right)_{x}= A
\left ( \begin{array}{ccc}
{\bf e}_{1} \\
{\bf e}_{2} \\
{\bf e}_{3}
\end{array} \right), \quad
\left ( \begin{array}{ccc}
{\bf e}_{1} \\
{\bf e}_{2} \\
{\bf e}_{3}
\end{array} \right)_{y}= B
\left ( \begin{array}{ccc}
{\bf e}_{1} \\
{\bf e}_{2} \\
{\bf e}_{3}
\end{array} \right) \eqno(11)
$$
where
$$
A =
\left ( \begin{array}{ccc}
0             & k     &  -\sigma \\
-\beta k      & 0     & \tau  \\
\sigma        & -\tau & 0
\end{array} \right) ,\quad
B =
\left ( \begin{array}{ccc}
0       & m_{3}  & -m_{2} \\
-\beta  m_{3} & 0      & m_{1} \\
\beta m_{2}  & -m_{1} & 0
\end{array} \right). \eqno(12)
$$

In fact, this equation is equivalent to the GMCE (10) with
$$
k=
\frac{L}{\sqrt{E}}, \quad \sigma =\frac{g}{E}\Gamma^{2}_{11},
\quad \tau = -\frac{g}{\sqrt{E}}p_{12}  \eqno(13a)
$$
$$
m_{1}= -\frac{g}{\sqrt{E}}p_{22},\quad
m_{2}= \frac{g}{E}\Gamma^{2}_{12}, \quad
m_{3}=
\frac{M}{\sqrt{E}} \eqno(13b)
$$
so that the compatibility condition of the equations (11) is the GMCE (10).
In this paper we usually suppose that the variables $k, \tau, \sigma,
m_{j}$ are some functions of $\lambda$ (the spectral parameter).

\section{Soliton geometry in $d=3$ dimensions}

\subsection{Curves and Surfaces in 2+1 dimensions}

It is known that the 2-dimensional SFE (11) or the GWE (7) admits
several (2+1)-dimensional integrable and non integrable generalizations,
which describe curves and surfaces in 2+1 dimensions.
Here some of them.

\subsubsection{The M-LIX equation}

The M-LIX  equation has the form [16]
$$
\alpha {\bf e}_{1y}=f_{1}{\bf e}_{1x}+
\sum_{j=1}^{n}b_{j}{\bf e}_{1}\wedge \frac{\partial^{j}}{\partial x^{j}}{\bf e}_{1} +
c_{1}{\bf e}_{2}+d_{1}{\bf e}_{3}  \eqno(14)
$$
and the  equations for ${\bf e}_{2y}, {\bf e}_{3y}, {\bf e}_{jt}$.
The M-LIX equation (14) admits several integrable reductions (see, e.g.
the subsection 3.3).

\subsubsection{The M-LXVII equation}
The M-LXVII equation looks like [16]
$$
\left ( \begin{array}{ccc}
{\bf e}_{1} \\
{\bf e}_{2} \\
{\bf e}_{3}
\end{array} \right)_{x}= A(\lambda)
\left ( \begin{array}{ccc}
{\bf e}_{1} \\
{\bf e}_{2} \\
{\bf e}_{3}
\end{array} \right) \eqno(15a)
$$
$$
\left ( \begin{array}{ccc}
{\bf e}_{1} \\
{\bf e}_{2} \\
{\bf e}_{3}
\end{array} \right)_{t} =  B(\lambda)
\left ( \begin{array}{ccc}
{\bf e}_{1} \\
{\bf e}_{2} \\
{\bf e}_{3}
\end{array} \right)_{y}+C(\lambda)
\left ( \begin{array}{ccc}
{\bf e}_{1} \\
{\bf e}_{2} \\
{\bf e}_{3}
\end{array} \right).
\eqno(15b)
$$
Here $A(\lambda), B(\lambda), C(\lambda)$ are some matrices.
This equation also contents some integrable reductions (see, e.g.
the subsection 3.5).

\subsubsection{The M-LX equation}

This extension of the GWE (7) or the SFE (1) has the form [16]
$$
\alpha\left ( \begin{array}{ccc}
{\bf e}_{1} \\
{\bf e}_{2} \\
{\bf e}_{3}
\end{array} \right)_{y}= A
\left ( \begin{array}{ccc}
{\bf e}_{1} \\
{\bf e}_{2} \\
{\bf e}_{3}
\end{array} \right)_{x}  + B
\left ( \begin{array}{ccc}
{\bf e}_{1} \\
{\bf e}_{2} \\
{\bf e}_{3}
\end{array} \right)
\eqno(16a)
$$
$$
\left ( \begin{array}{ccc}
{\bf e}_{1} \\
{\bf e}_{2} \\
{\bf e}_{3}
\end{array} \right)_{t}=
\sum_{j=0}^{n}C_{j}\frac{\partial^{j}}{\partial x^{j}}
\left ( \begin{array}{ccc}
{\bf e}_{1} \\
{\bf e}_{2} \\
{\bf e}_{3}
\end{array} \right)
\eqno(16b)
$$
where $A, B, C_{j}$ - some matrices. Some integrable reductions of the
M-LX equation (16) were presented in [16].

\subsubsection{The modified M-LXI equation}

The modified M-LXI (mM-LXI) equation, which is the (2+1)-dimensional
integrable extension of the GWE (7), usually we write in the form [16]
$$
\left ( \begin{array}{c}
{\bf e}_{1} \\
{\bf e}_{2} \\
{\bf e}_{3}
\end{array} \right)_{x}= A
\left ( \begin{array}{ccc}
{\bf e}_{1} \\
{\bf e}_{2} \\
{\bf e}_{3}
\end{array} \right), \quad
\left ( \begin{array}{ccc}
{\bf e}_{1} \\
{\bf e}_{2} \\
{\bf e}_{3}
\end{array} \right)_{y}= B
\left ( \begin{array}{ccc}
{\bf e}_{1} \\
{\bf e}_{2} \\
{\bf e}_{3}
\end{array} \right), \quad
\left ( \begin{array}{ccc}
{\bf e}_{1} \\
{\bf e}_{2} \\
{\bf e}_{3}
\end{array} \right)_{t}= C
\left ( \begin{array}{ccc}
{\bf e}_{1} \\
{\bf e}_{2} \\
{\bf e}_{3}
\end{array} \right)
\eqno(17)
$$
with
$$
A =
\left ( \begin{array}{ccc}
0             & k     &  -\sigma \\
-\beta k      & 0     & \tau  \\
\sigma        & -\tau & 0
\end{array} \right) ,
B=
\left ( \begin{array}{ccc}
0            & m_{3}  & -m_{2} \\
-\beta m_{3} & 0      & m_{1} \\
\beta m_{2}  & -m_{1} & 0
\end{array} \right),
C =
\left ( \begin{array}{ccc}
0       & \omega_{3}  & -\omega_{2} \\
-\beta\omega_{3} & 0      & \omega_{1} \\
\beta\omega_{2}  & -\omega_{1} & 0
\end{array} \right). \eqno(18)
$$
In the case  $\sigma = 0$ the mM-LXI equation transform to the M-LXI
equation. We think that the M-LXI equation (17) is integrable in general
and admits integrable reductions in particular (see, e.g.  the
subsection 3.2).

\subsection{The mM-LXII equation and Soliton equations in 2+1 dimensions}

Let us return to the mM-LXI equation (17).
From (17) we obtain the following  mM-LXII equation [16]
$$
A_{y} - B_{x} + [A,B] = 0 \eqno(19a)
$$
$$
A_{t}-C_{x}+[A,C]=0   \eqno(19b)
$$
$$
B_{t}-C_{y}+[B,C]=0. \eqno(19c)
$$

The mM-LXI and mM-LXII equations are may be some of the
 (2+1)-dimensional simplest, interesting and important integrable
extensions of the GWE (7)=(11) and GMCE (10), respectively.
In fact, we observe that  the M-LXII equation (19) is
the particular case of the other
integrable equation, namely the Bogomolny equation (BE). The BE has the form
(see, e.g. [15])
$$
\Phi_{t}+[\Phi,C]+ A_{y} - B_{x} + [A,B] = 0 \eqno(20a)
$$
$$
\Phi_{y}+[\Phi,B]+C_{x}-A_{t}+[C,A]=0   \eqno(20b)
$$
$$
\Phi_{x}+[\Phi,A]+B_{t}-C_{y}+[B,C]=0. \eqno(20c)
$$
Hence as $\Phi=0$ we obtain the mM-LXII (or M-LXII) equation (19). Hence
follows and the other remarkable fact, namely, the mM-LXII equation
is exact reduction of the Self-Dual Yang-Mills equation (SDYME)
$$
F_{\alpha\beta}=0, \quad F_{\bar\alpha\bar\beta}=0,
\quad
F_{\alpha\bar\alpha}+F_{\beta\bar\beta}=0  \eqno(21)
$$
Here
$$
F_{\mu\nu}=\frac{\partial A_{\nu}}{\partial x_{\mu}}-
\frac{\partial A_{\mu}}{\partial x_{\nu}}+[A_{\mu},A_{\nu}] \eqno(22)
$$
and
$$
\frac{\partial}{\partial x_{\alpha}}=\frac{\partial}{\partial z}-
i\frac{\partial}{ \partial t}, \quad
\frac{\partial}{\partial x_{\bar\alpha}}=\frac{\partial}{\partial z}+
i\frac{\partial}{ \partial t}, \quad
\frac{\partial}{\partial x_{\beta}}=\frac{\partial}{\partial x}-
i\frac{\partial}{ \partial y}, \quad
\frac{\partial}{\partial x_{\bar\beta}}=\frac{\partial}{\partial x}+
i\frac{\partial}{ \partial y}. \eqno(23)
$$
If in the SDYME (21) we take
$$
A_{\alpha}= -iC, \quad A_{\bar\alpha}= iC, \quad
A_{\beta}= A-iB, \quad A_{\beta}= A+iB \eqno(100)
$$
and if $A,B,C$ are independent of $z$, then the SDYME (21)
reduces to  the mM-LXII (or M-LXII as $\sigma=0$) equation (19).
As known that the LR of the SDYME has the form [28, 29]
$$
(\partial_{\alpha}+\lambda\partial_{\bar\beta})\Psi=(A_{\alpha}+\lambda A_{\bar\beta})\Psi,
\quad
(\partial_{\beta}-\lambda\partial_{\bar\alpha})\Psi=(A_{\beta}-
\lambda A_{\bar\alpha})\Psi \eqno(24)
$$
where $\lambda$ is the spectral parameter  satisfing the following
set of the equations
$$
\lambda_{\beta}=\lambda\lambda_{\bar\alpha}, \quad
\lambda_{\alpha}=-\lambda\lambda_{\bar\beta}. \eqno(25)
$$
Apropos, the simplest solution of this set has may be the following
form [16]
$$
\lambda =\frac{a_{1}x_{\bar\alpha}+a_{2}x_{\bar\beta}+a_{3}}{a_{2}x_{\alpha}-
a_{1}x_{\beta}+a_{4}}, \quad a_{j}=consts.
$$
From (24) we obtain the LR of the M-LXII equation (19)
$$
(-i\partial_{t}+\lambda\partial_{\bar\beta})\Psi=[-iC+\lambda (A+iB)]\Psi
\eqno(26a)
$$
$$
(\partial_{\beta}-i\lambda\partial_{t})\Psi=[(A-iB)-i\lambda C]\Psi.
 \eqno(26b)
$$
So the mM-LXII equation (19) is integrable in the sense that it admits the LR (26).
The other form LR for it we present in subsection 3.4 (see (62)).

Now we will establish the connection between the
M-LXI and mM-LXII equations (17), (19)  and soliton equations in 2+1 dimensions.
Let us, we assume
$$
{\bf e}_{1} \equiv {\bf S}. \eqno(27)
$$
Moreover we introduce two complex functions
$q, p$ according to the following expressions
$$
q = a_{1}e^{ib_{1}}, \quad p=a_{2}e^{ib_{2}}  \eqno(28)
$$
where $a_{j}, b_{j}$ are real functions. Now we ready to consider some
examples.

\subsubsection{ The Ishimori and Davey-Stewartson equations}.

The Ishimori equation (IE) reads as
$$
{\bf S}_{t}  =  {\bf S}\wedge ({\bf S}_{xx} +\alpha^2 {\bf S}_{yy})+
u_x{\bf S}_{y}+u_y{\bf S}_{x} \eqno (29a)
$$
$$
u_{xx}-\alpha^2 u_{yy}  = -2\alpha^2 {\bf S}\cdot ({\bf S}_{x}\wedge
                        {\bf S}_{y}).   \eqno (29b)
$$
In this case we have [21]
$$
m_{1}=\partial_{x}^{-1}[\tau_{y}-\frac{\epsilon}{2\alpha^2}M_2^{Ish}u],
\quad
m_{2}= -\frac{1}{2\alpha^2 k}M_2^{Ish}u, \quad
m_{3} =\partial_{x}^{-1}[k_y +\frac{\tau}{2\alpha^2 k}M_2^{Ish}u] \eqno(30)
$$
and
$$
\omega_{1} = \frac{1}{k}[-\omega_{2x}+\tau\omega_{3}],\quad
\omega_{2}= -k_{x}-
 \alpha^{2}(m_{3y}+m_{2}m_{1})+im_{2}u_{x}
$$
$$
\omega_{3}= -k \tau+\alpha^{2}(m_{2y}-m_{3}m_{1})
+ik u_{y}+im_{3}u_{x}, \quad
M_2^{Ish}=M_2|_{a=b=-\frac{1}{2}}. \eqno(31)
$$
Functions $q, p$ are given by (28) with
$$
a_{1}^2 =a_{1}^{\prime^{2}}=\frac{1}{4}k^2+
\frac{|\alpha|^2}{4}(m_3^2 +m_2^2)-\frac{1}{2}\alpha_{R}km_3-
\frac{1}{2}\alpha_{I}km_2
  \eqno(32a)
$$
$$
b_1 =\partial_{x}^{-1}\{-\frac{\gamma_1}{2ia_1^{\prime^{2}}}-(\bar A-A+D-\bar D)\}  \eqno(32b)
$$
$$
a_2^2=a_{2}^{\prime^{2}}=\frac{1}{4}k^2+
\frac{|\alpha|^2}{4}(m_3^2 +m_2^2)+\frac{1}{2}\alpha_{R}km_3
-\frac{1}{2}\alpha_{I}km_2
  \eqno(32c)
$$
$$
b_{2} =\partial_{x}^{-1}\{-\frac{\gamma_2}{2ia_2^{\prime^{2}}}-(A-\bar A+\bar D-D)\}  \eqno(32d)
$$
where
$$
\gamma_1=i\{\frac{1}{2}k^{2}\tau+
\frac{|\alpha|^2}{2}(m_3km_1+m_2k_y)-
$$
$$
\frac{1}{2}\alpha_{R}(k^{2}m_1+m_3k\tau+
m_2k_x)
+\frac{1}{2}\alpha_{I}[k(2k_y-m_{3x})-
k_x m_3]\} \eqno(33a)
$$
$$
\gamma_2=-i\{\frac{1}{2}k^{2}\tau+
\frac{|\alpha|^2}{2}(m_3km_1+m_2 k_y)+
$$
$$
\frac{1}{2}\alpha_{R}(k^{2}m_1+m_3k\tau+
m_2k_x )
+\frac{1}{2}\alpha_{I}[k(2k_y-m_{3x})-
k_x m_3]\}. \eqno(33b)
$$
Here $\alpha=\alpha_{R}+i\alpha_{I}$. In this case, $q,p$ satisfy the following
DS  equation
$$
iq_t  + q_{xx}+\alpha^{2}q_{yy} + vq = 0 \eqno (34a)
$$
$$
-ip_t +  p_{xx}+\alpha^{2}p_{yy} + vp = 0 \eqno (34b)
$$
$$
v_{xx}-\alpha^{2}v_{yy} + 2[(p q)_{xx}+\alpha^{2}(p q)_{yy}] = 0.
\eqno (34c)
$$
So we have proved that the IE (29) and the DS equation (34)
are L-equivalent to each other.
As well known that these equations are G-equivalent to each other [31].
Note that the IE contains two reductions:
the Ishimori I equation as $\alpha_{R}=1, \alpha_{I}=0$ and
the Ishimori II equation as $\alpha_{R}=0, \alpha_{I}=1$.
The corresponding versions of the DS equation (34), we obtain as
the corresponding values of the parameter $\alpha$ (for details, see, e.g.
the ref. [21]).

\subsubsection{ The Zakharov II and  M-IX  equations}

Now we find the connection between the Myrzakulov IX (M-IX) equation
and the curves (the M-LXI equation). The M-IX equation reads as
$$
{\bf S}_t = {\bf S} \wedge M_1{\bf S}+A_2{\bf S}_x+A_1{\bf S}_y  \eqno(35a)
$$
$$
M_2u=2\alpha^{2} {\bf S}({\bf S}_x \wedge {\bf S}_y) \eqno(35b)
$$
where $ \alpha,b,a  $=  consts and
$$
M_1= \alpha ^2\frac{\partial ^2}{\partial y^2}+4\alpha (b-a)\frac{\partial^2}
   {\partial x \partial y}+4(a^2-2ab-b)\frac{\partial^2}{\partial x^2}
$$
$$
M_2=\alpha^2\frac{\partial^2}{\partial y^2} -2\alpha(2a+1)\frac{\partial^2}
   {\partial x \partial y}+4a(a+1)\frac{\partial^2}{\partial x^2}
$$
$$
A_1=i\{\alpha (2b+1)u_y - 2(2ab+a+b)u_{x}\}
$$
$$
A_2=i\{4\alpha^{-1}(2a^2b+a^2+2ab+b)u_x - 2(2ab+a+b)u_{y}\}. \eqno(36)
$$
The M-IX equation was introduced in [16] and is integrable. It admits several
integrable reductions:

i) the Ishimori equation (29) as $a=b=-\frac{1}{2}$;

ii) the M-VIII equation as $a=b=-1, X=x/2, Y=y/\alpha, w=-\alpha^{-1}u_{Y}$
$$
{\bf S}_t = {\bf S} \wedge {\bf S}_{YY}+w{\bf S}_{Y}  \eqno(35a)
$$
$$
w_{X}+w_{Y} + {\bf S}({\bf S}_{X} \wedge {\bf S}_{Y})=0 \eqno(35b)
$$

iii) the M-XXXIV equation as $a=b=-1, X=t$
$$
{\bf S}_t = {\bf S} \wedge {\bf S}_{YY}+w{\bf S}_{Y}  \eqno(35a)
$$
$$
w_{t}+w_{Y} + \frac{1}{2}({\bf S}^{2}_{Y})_{Y}=0 \eqno(35b)
$$
and so on [16].
In our  case we have
$$
m_{1}=\partial_{x}^{-1}[\tau_{y}-\frac{\beta}{2\alpha^2}M_2 u],\quad
m_{2}=-\frac{1}{2\alpha^2 k}M_2 u, \quad
m_{3}=\partial_{x}^{-1}[k_y +\frac{\tau}{2\alpha^2 k}M_2 u]  \eqno(37)
$$
and
$$
\omega_{1} = \frac{1}{k}[-\omega_{2x}+\tau\omega_{3}],
                 \eqno (38a)
$$
$$
\omega_{2}= -4(a^{2}-2ab-b)k_{x}-
4\alpha (b-a)k_{y} -\alpha^{2}(m_{3y}+m_{2}m_{1})+m_{2}A_{1}
\eqno (38b)
$$
$$
\omega_{3}= -4(a^{2}-2ab-b)k \tau-
4\alpha (b-a)k m_{1}+\alpha^{2}(m_{2y}-m_{3}m_{1})
+k A_{2}+m_{3}A_{1}.
\eqno (38c)
$$
Functions $q, p$ are given by (28) with
$$
a_{1}^2 =\frac{|a|^2}{|b|^2}a_{1}^{\prime^{2}}=\frac{|a|^2}{|b|^2}\{(l+1)^2k^2
+\frac{|\alpha|^2}{4}(m_3^2 +m_2^2)-(l+1)\alpha_{R}km_3-
(l+1)\alpha_{I}km_2\}
  \eqno(39a)
$$
$$
b_{1} =\partial_{x}^{-1}\{-\frac{\gamma_1}{2ia_1^{\prime^{2}}}-(\bar A-A+D-\bar D)\}  \eqno(39b)
$$
$$
a_{2}^2 =\frac{|b|^2}{|a|^2}a_{2}^{\prime^{2}}=\frac{|b|^2}{|a|^2}\{l^2k^2
+\frac{|\alpha|^2}{4}(m_3^2 +m_2^2)-l\alpha_{R}km_3+
l\alpha_{I}km_2\}
  \eqno(39c)
$$
$$
b_{2} =\partial_{x}^{-1}\{-\frac{\gamma_2}{2ia_2^{\prime^{2}}}-(A-\bar A+\bar D-D)  \eqno(39d)
$$
where
$$
\gamma_1=i\{2(l+1)^2k^{2}\tau+\frac{|\alpha|^2}{2}(m_3km_1+m_2k_y)-
$$
$$
(l+1)\alpha_{R}[k^{2}m_1+m_3k\tau+
m_2k_x]+(l+1)\alpha_{I}[k(2k_y-m_{3x})-
k_x m_3]\} \eqno(40a)
$$
$$
\gamma_2=-i\{2l^2k^{2}\tau+
\frac{|\alpha|^2}{2}(m_3km_1+m_2k_y)-
$$
$$
l\alpha_{R}(k^{2}m_1+m_3k\tau+
m_2k_x)-l\alpha_{I}[k(2k_y-m_{3x})-
k_x m_3]\}. \eqno(40b)
$$
Here $\alpha=\alpha_{R}+i\alpha_{I}$. In this case, $q,p$
satisfy the following Zakharov II (Z-II) equation [32]
$$
iq_t+M_{1}q+vq=0 \eqno(41a)
$$
$$
ip_t-M_{1}p-vp=0 \eqno(41b)
$$
$$
M_{2}v=-2M_{1}(pq). \eqno(41c)
$$

As well known the Z-II equation admits several reductions:
1)  the DS-I equation as $\alpha_{R}=1, \alpha_{I}=0$;
2)  the DS-II equation as $\alpha_{R}=0, \alpha_{I}=1$;
3) the Z-III equation as $a=b=-1$
and so on [32, 16].

\subsection{The M-LIX equation and Soliton equations in 2+1 dimensions}

Now let us consider the connection between the M-LIX equation (14) and
(2+1)-dimensional soliton equations. Mention that the M-LIX equation is
one of (2+1)-dimensional extensions of the SFE (1).
As example, let us consider the connection between the M-LIX equation
and the Z-II equation (41). Let the spatial part of the M-LIX equation
has the form [16]
$$
\alpha {\bf e}_{1\eta}=i{\bf e}_{1}\wedge{\bf e}_{1\xi} +
+i(q+p){\bf e}_{2}+(q-p){\bf e}_{3}  \eqno(42a)
$$
and the equations for the
$ {\bf e}_{2y}, \quad  {\bf e}_{3y}.$ This set of these equations can
be considered as some generalization of the
Belavin-Polyakov  equation [28]
$$
{\bf e}_{1\eta}=\pm {\bf e}_{1}\wedge{\bf e}_{1\xi}
  \eqno(43)
$$
which arises in several physical applications.
In terms of matrices the equations (42)  we
can write in the form [16]
$$
\alpha \hat  e_{1\eta} = \frac{1}{2}[\hat e_{1},\hat e_{1\xi}] +
i(q+p)\hat e_{2}+(q-p)\hat e_{3}  \eqno(44a)
$$
$$
\alpha \hat  e_{2\eta} =  [\hat e_{1}, \hat e_{2\xi}]+
i\hat e_{3\xi}+i(A+B)\hat e_{3}+
(A-B)\hat e_{2}-i(p+q)\hat e_{1}  \eqno(44b)
$$
$$
\alpha \hat  e_{3\eta} = [\hat e_{1}, \hat e_{3\xi}]-i\hat e_{2\xi}-i(A+B)\hat e_{2}+
(A-B)\hat e_{3}+(p-q)\hat e_{1}   \eqno(44c)
$$
where
$$
\hat  e_{1} = g^{-1}\sigma_{3}g, \quad \hat e_{2}=g^{-1}\sigma_{2}g,
\quad \hat e_{3} = g^{-1}\sigma_{1}g, \quad \xi =\frac{y}{\alpha},
\quad \eta =2x + \frac{2a+1}{\alpha}y.  \eqno(45)
$$
Here  $\sigma_{j}$ are Pauli matrices.
Hence follows that the matrix-function  $g$ satisfies the equations
$$
\alpha g_{y}=  B_{1}g_{x} + B_{0}g \eqno(46)
$$
with
$$
B_{0}=
\left( \begin{array}{cc}
0  &  q\\
p  &  0
\end{array} \right), \quad
B_{1}=
\left ( \begin{array}{cc}
a+1  &  0\\
0    &  a
\end{array} \right ).  \eqno(47)
$$
To find the time evolution of matrices  $\hat e_{j}$
we require that the matrices  $\hat e_{j}$
satisfy the  following set of the equations which is the time part of the
M-LIX equation (14)
$$
i\hat  e_{1t} =-\{[(2b+1)\hat e_{1}+1][\hat e_{1},\hat e_{1\xi\xi}]+
[2(2b+1)F^{+}+F^{-}][\hat e_{1},\hat e_{1\xi}]+Q\hat e_{1\xi}\}+
$$
$$
-\{[(c_{12}+c_{21})+i(p-q)F^{+}]\hat e_{2}-[i(c_{12}-c_{21})+(p-q)F^{-}]\hat e_{3}\}
\eqno(48a)
$$
$$
i\hat  e_{2t} =[-P+i(c_{11}-c_{22})]\hat e_{3}+
[i(p-q)F^{+}+(c_{12}+c_{21})]\hat e_{1}+
-Q \hat e_{3}\hat e_{1\xi} +i\hat e_{3}\hat e_{1\xi\xi} \hat e_{1} +
i(p-q)\hat e_{1}\hat e_{1\xi}+
$$
$$
\frac{i}{4}[i(2b+1)\hat e_{2}-\hat e_{3}][\hat e_{1},\hat e_{1\xi\xi}]-
T[\hat e_{2},\hat e_{1\xi}]+[2b+1+\hat e_{1}]\{i\hat e_{1\xi\xi}\hat e_{3}+
\frac{1}{4}[\hat e_{1},\hat e_{1\xi\xi}]\hat e_{2}+
\frac{1}{2}[\hat e_{2},\hat e_{1\xi\xi}]\hat e_{1}\}
\eqno(48b)
$$
$$
i\hat  e_{3t} = P\hat e_{2}+[c_{11}-c_{22}-i(c_{12}-c_{21})-(p-q)F^{-}]\hat e_{1}
+ Q\hat e_{2}\hat e_{1\xi}+T[\hat e_{3},\hat e_{1\xi}]-
\frac{1}{4}[(2b+1)\hat e_{3}-i\hat e_{2}][\hat e_{1},\hat e_{1\xi\xi}-
$$
$$
i\hat e_{2}\hat e_{1\xi\xi}\hat e_{1}-[2b+1+\hat e_{1}]\{i\hat e_{1\xi\xi}\hat e_{2}
+\frac{1}{4}[\hat e_{1},\hat e_{1\xi\xi}]\hat e_{3}+
\frac{1}{2}[\hat e_{3},\hat e_{1\xi\xi}]\hat e_{1}\}
\eqno(48c)
$$
where
$$
Q=2(2b+1)F^{-}+4F^{+}+(p+q),\quad P=i[2(2b+1)(F^{-}F^{+}+F^{+}_{\xi})+
(p+q)F^{+}+F^{{-}^{2}}+F^{{+}^{2}}+2F^{-}_{\xi}],
$$
$$
T=2(2b+1)F^{+}+(2b+1)F^{-}\hat e_{1}+2F^{+}+F^{-}+\frac{1}{2}(p+q)\hat e_{1}+
\frac{1}{2}(p-q)\hat e_{3}.
$$
We note that as  follows from these equations the vector ${\bf e}_{1}$
satisfies the M-IX equation
$$
i\hat  e_{1t} = \frac{1}{2}[\hat e_{1}, M_{1}\hat e_{1}] +A_{1}\hat e_{1y}+A_{2}\hat e_{1x}
\eqno(49a)
$$
$$
M_{2}u = \frac{\alpha^{2}}{2i}tr(\hat e_{1}([\hat e_{1x},\hat e_{1y}]).
\eqno(49b)
$$

The time part of the M-LIX equation (48) immediently gives
the equation
$$
g_{t}=  2C_{2}g_{xx} + C_{1}g_{x} +C_{0}g \eqno(50)
$$
where
$$
C_{0} =
\left ( \begin{array}{cc}
c_{11}      &  c_{12} \\
c_{21}       &  c_{22}
\end{array} \right), \quad
C_{2}=\frac{2b+1}{2}I+\frac{1}{2}\sigma_{3},
\quad C_{1}=iB_{0}.    \eqno(51)
$$
Here
$$
c_{12}=i(2b-a+1)q_{x}+i\alpha q_{y},\quad c_{21}=i(a-2b)q_{x}-i\alpha p_{y} \eqno(52)
$$
and $c_{jj}$ are the solutions of the following equations
$$
(a+1)c_{11x}-\alpha c_{11y}=i[(2b-a+1)(pq)_{x}+\alpha(pq)_{y}]  \eqno(53a)
$$
$$
ac_{22x}-\alpha c_{22y}=i[(a-2b)(pq)_{x}-\alpha (pq)_{y}].       \eqno(53b)
$$

So we have identified the curve,
given by the M-LIX equation (42) and (48) with
the M-IX equation (35)$\equiv$(49). On the other hand, the compatibilty
condition of equations (46) and (50) is equivalent to the Z-II equation
(41).
So that we have also established the connection between the  curve
(the M-LIX equation) and the Z-II equation.
And we have shown, once more that
the M-IX equation (35) and the Z-II equation (41) are L-equivalent to each other.
Finally we note as $a=b=-\frac{1}{2}$ from these results follow
the corresponding connection between the M-LIX, Ishimori and DS equations [21].
And as $a=b=-1$ we get the relation between the
M-VIII, M-LIX and Zakharov-III equations (for details, see [16]).

\subsection{The M-LVIII, M-LXIII  equations and Soliton equationa in 2+1 dimensions}

In the C-approach, our starting point is  the following (2+1)-dimensional
M-LVIII equation [16]
$$
{\bf r}_{t} = \Upsilon_{1} {\bf r}_{x} + \Upsilon_{2} {\bf r}_{y}
+ \Upsilon_{3}{\bf n}   \eqno(54a)
$$
$$
{\bf r}_{xx} = \Gamma^{1}_{11} {\bf r}_{x} + \Gamma^{2}_{11} {\bf r}_{y}
+ L{\bf n}   \eqno(54b)
$$
$$
{\bf r}_{xy} = \Gamma^{1}_{12} {\bf r}_{x} + \Gamma^{2}_{12} {\bf r}_{y}
+ M{\bf n}   \eqno(54c)
$$
$$
{\bf r}_{yy} = \Gamma^{1}_{22} {\bf r}_{x} + \Gamma^{2}_{22} {\bf r}_{y}
+ N{\bf n}   \eqno(54d)
$$
$$
{\bf n}_{x} = p_{11} {\bf r}_{x} + p_{12} {\bf r}_{y}   \eqno(54e)
$$
$$
{\bf n}_{y} = p_{21} {\bf r}_{x} + p_{22} {\bf r}_{y}.   \eqno(54f)
$$
It is one of the (2+1)-dimensional generalizations of the GWE (3)
and admits  several integrable reductions. Practically, all
integrable spin systems  in 2+1 dimensions are some
integrable reductions
of the M-LVIII equation (54). For eexample, the isotropic
M-I equation (76) or in the vector form
$$
{\bf S}_t = ({\bf S} \wedge {\bf S}_{y}+u{\bf S})_x  \eqno(101a)
$$
$$
u_{x}=-{\bf S}({\bf S}_x \wedge {\bf S}_y) \eqno(101b)
$$
is the particular case of the M-LVIII equation (54). In fact,
let ${\bf r}_{x}={\bf S}$, then the M-I equation (101) takes the form
$$
{\bf r}_{t} = (u+\frac{MF}{\sqrt{g}}) {\bf r}_{x} -
\frac{M}{\sqrt{g}} {\bf r}_{y} + \Gamma^{2}_{12} {\bf n}   \eqno(102)
$$
with
$$
u=\partial^{-1}_{x}[\sqrt{g}(L\Gamma^{2}_{12}-M\Gamma^{2}_{11}). \eqno(103)
$$
This equation (102) is in fact the particular case of the M-LVIII
equation (54) with
$$
\Upsilon_{1}=u+\frac{MF}{\sqrt{g}}, \quad
\Upsilon_{2}= -\frac{M}{\sqrt{g}},\quad
\Upsilon_{3}= \Gamma^{2}_{12} \sqrt{g}.   \eqno(104)
$$

Sometimes it is convenient to work using the B-approach. In this approach
the starting equation is the  following M-LXIII equation [16]
$$
{\bf r}_{tx} = \Gamma^{1}_{01} {\bf r}_{x} + \Gamma^{2}_{01} {\bf r}_{y}
+ \Gamma^{3}_{01}{\bf n}   \eqno(55a)
$$
$$
{\bf r}_{ty} = \Gamma^{1}_{02} {\bf r}_{x} + \Gamma^{2}_{02} {\bf r}_{y}
+ \Gamma^{3}_{02}{\bf n}   \eqno(55b)
$$
$$
{\bf r}_{xx} = \Gamma^{1}_{11} {\bf r}_{x} + \Gamma^{2}_{11} {\bf r}_{y}
+ L{\bf n}   \eqno(55c)
$$
$$
{\bf r}_{xy} = \Gamma^{1}_{12} {\bf r}_{x} + \Gamma^{2}_{12} {\bf r}_{y}
+ M{\bf n}   \eqno(55d)
$$
$$
{\bf r}_{yy} = \Gamma^{1}_{22} {\bf r}_{x} + \Gamma^{2}_{22} {\bf r}_{y}
+ N{\bf n}   \eqno(55e)
$$
$$
{\bf n}_{t} = p_{01} {\bf r}_{x} + p_{02} {\bf r}_{y}   \eqno(55f)
$$
$$
{\bf n}_{x} = p_{11} {\bf r}_{x} + p_{12} {\bf r}_{y}   \eqno(55g)
$$
$$
{\bf n}_{y} = p_{21} {\bf r}_{x} + p_{22} {\bf r}_{y}.   \eqno(55h)
$$
This equation follows from the M-LVIII equation (54) under the following conditions
$$
\Gamma^{1}_{01}= \Upsilon_{1x}+ \Upsilon_{1} \Gamma^{1}_{11}
+ \Upsilon_{2}\Gamma^{1}_{12}+\Upsilon_{3}p_{11}
$$
$$
\Gamma^{1}_{02}= \Upsilon_{2x}+ \Upsilon_{1} \Gamma^{2}_{11}
+ \Upsilon_{2}\Gamma^{2}_{12}+\Upsilon_{3}p_{12}
$$
$$
\Gamma^{1}_{03}= \Upsilon_{3x}+ \Upsilon_{1} L
+ \Upsilon_{2}M
$$
$$
p_{01}= \frac{F \Gamma^{3}_{02}}{g}, \quad
p_{02}= -\frac{E \Gamma^{3}_{02}}{g}, \quad
g = EG-F^{2} \eqno(56)
$$
Note that the M-LXIII equation (69) usually we use in the following form
$$
Z_{x} = AZ, \quad Z_{y} = BZ, \quad Z_{t} = CZ  \eqno(57)
$$
where
$$
A =
\left ( \begin{array}{ccc}
\Gamma^{1}_{11} & \Gamma^{2}_{11} & L \\
\Gamma^{1}_{12} & \Gamma^{2}_{12} & M \\
p_{11}           & p_{12}           & 0
\end{array} \right) ,  \quad
B =
\left ( \begin{array}{ccc}
\Gamma^{1}_{12} & \Gamma^{2}_{12} & M \\
\Gamma^{1}_{22} & \Gamma^{2}_{22} & N \\
p_{21}          & p_{22}          & 0
\end{array} \right), \quad
C =
\left ( \begin{array}{ccc}
\Gamma^{1}_{01} & \Gamma^{2}_{01} & \Gamma^{3}_{01} \\
\Gamma^{1}_{02} & \Gamma^{2}_{02} & \Gamma^{3}_{02} \\
\Gamma^{1}_{03} & \Gamma^{2}_{03} & 0
\end{array} \right) . \eqno(58)
$$

The compatibility condition of the equations (57) gives
$$
A_{y}-B_{x}+[A,B]=0 \eqno(59a)
$$
$$
A_{t}-C_{x}+[A,C]=0 \eqno(59b)
$$
$$
B_{t}-C_{y}+[B,C]=0 \eqno(59c)
$$
that is  the mM-LXII equation (19). These equations are equivalent the relations
$$
{\bf r}_{yxx} = {\bf r}_{xxy}, \quad
{\bf r}_{yyx}={\bf r}_{xyy}  \eqno(60a)
$$
$$
{\bf r}_{txx} = {\bf r}_{xxt}, \quad
{\bf r}_{txy}={\bf r}_{xyt}, \quad {\bf r}_{tyy}={\bf r}_{yyt}.  \eqno(60b)
$$
Note that (60a) is the well known GMCE (10). So that the M-LXII equation
(59) is one of the (2+1)-dimensional generalizations of the GMCE (10).
In the  orthogonal basis (6) the equation (57) takes the form
$$
\left ( \begin{array}{ccc}
{\bf e}_{1} \\
{\bf e}_{2} \\
{\bf e}_{3}
\end{array} \right)_{x}=\frac{1}{\sqrt{E}}
\left ( \begin{array}{ccc}
0             & L   & -\frac{g}{\sqrt{E}}\Gamma^{2}_{11} \\
-\beta L      & 0     & -g p_{12}  \\
\frac{\beta g}{\sqrt{E}}\Gamma^{2}_{11}& g p_{12} & 0
\end{array} \right)
\left ( \begin{array}{c}
{\bf e}_{1} \\
{\bf e}_{2} \\
{\bf e}_{3}
\end{array} \right)
\eqno(61a)
$$
$$
\left ( \begin{array}{ccc}
{\bf e}_{1} \\
{\bf e}_{2} \\
{\bf e}_{3}
\end{array} \right)_{y}=\frac{1}{\sqrt{E}}
\left ( \begin{array}{ccc}
0             & M   & -\frac{g}{\sqrt{E}}\Gamma^{2}_{12} \\
-\beta M      & 0     & -g p_{22}  \\
\frac{\beta g}{\sqrt{E}}\Gamma^{2}_{12}& g p_{22} & 0
\end{array} \right)
\left ( \begin{array}{c}
{\bf e}_{1} \\
{\bf e}_{2} \\
{\bf e}_{3}
\end{array} \right)
\eqno(61b)
$$
$$
\left ( \begin{array}{ccc}
{\bf e}_{1} \\
{\bf e}_{2} \\
{\bf e}_{3}
\end{array} \right)_{t}=\frac{1}{\sqrt{E}}
\left ( \begin{array}{ccc}
0             & \Gamma^{3}_{01}   & -\frac{g}{\sqrt{E}}\Gamma^{2}_{01} \\
-\beta\Gamma^{3}_{01} & 0     & -g \Gamma^{2}_{03}  \\
\frac{\beta g}{\sqrt{E}}\Gamma^{2}_{01}& g \Gamma^{2}_{03} & 0
\end{array} \right)
\left ( \begin{array}{c}
{\bf e}_{1} \\
{\bf e}_{2} \\
{\bf e}_{3}
\end{array} \right).
\eqno(61c)
$$

This equations we can rewrite in terms of 2$\times$2 matrices as
$$
g_{x}=Ug,\quad g_{y}=Vg, \quad  g_{t}=Wg   \eqno(62)
$$
where
$$
U =
\frac{1}{2i\sqrt{E}}
\left ( \begin{array}{cc}
-\sqrt{g} p_{12}  & L+i\sqrt{\frac{g}{E}}\Gamma^{2}_{11} \\
L-i\sqrt{\frac{g}{E}}\Gamma^{2}_{11} &\sqrt{g} p_{12}
\end{array} \right)
\eqno(63a)
$$
$$
V=
\frac{1}{2i\sqrt{E}}
\left ( \begin{array}{cc}
-\sqrt{g} p_{22}  & M-i\sqrt{\frac{g}{E}}\Gamma^{2}_{12} \\
M+i\sqrt{\frac{g}{E}}\Gamma^{2}_{12} &\sqrt{g} p_{22}
\end{array} \right)
\eqno(63b)
$$
$$
W=
\frac{1}{2i\sqrt{E}}
\left ( \begin{array}{cc}
-\sqrt{g}\Gamma^{2}_{03}  & \Gamma^{3}_{01}-i\sqrt{\frac{g}{E}}\Gamma^{2}_{01} \\
\Gamma^{3}_{01}+i\sqrt{\frac{g}{E}}\Gamma^{2}_{01} &\sqrt{g} \Gamma^{2}_{03}
\end{array} \right).
\eqno(63c)
$$

From these equations follow
$$
U_{y}-V_{x}+[U,V]=0 \eqno(64a)
$$
$$
U_{t}-W_{x}+[U,W]=0 \eqno(64b)
$$
$$
V_{t}-W_{y}+[V,W]=0 \eqno(64c)
$$
that  is the other form of the mM-LXII equation (59).
The equation (64a) is the GMCE (10).
Note that the M-LXIII equation in the form  (61) have
the same form with the mM-LXI equation (17) with the following
identifications
$$
k=
\frac{L}{\sqrt{E}}, \quad \sigma =\frac{g}{E}\Gamma^{2}_{11},
\quad \tau = -\frac{g}{\sqrt{E}}p_{12}  \eqno(65a)
$$
$$
m_{1}= -\frac{g}{\sqrt{E}}p_{22},\quad
m_{2}= \frac{g}{E}\Gamma^{2}_{12}, \quad
m_{3}=
\frac{M}{\sqrt{E}} \eqno(65b)
$$
$$
\omega_{1}=-\frac{1}{\sqrt{E}}g \Gamma^{2}_{03}, \quad
\omega_{2}=\frac{g}{E}\Gamma^{2}_{03},
\quad
\omega_{3}=\frac{1}{\sqrt{E}}\Gamma^{3}_{01}. \eqno(65c)
$$
So that the set of the linear equations (62) can be considered as one of the
form of the LR for the mM-LXII equation (19).

\subsection{The M-LXVII equation and Soliton equations in 2+1 dimensions}

In this section we consider curves and/or surfaces which are given by
the following M-LXVII equation [16]
$$
\left( \begin{array}{c}
{\bf e}_{1}\\
{\bf e}_{2}  \\
{\bf e}_{3}
\end{array} \right)_{\xi_1} = B
\left( \begin{array}{c}
{\bf e}_{1}\\
{\bf e}_{2}  \\
{\bf e}_{3}
\end{array} \right), \quad
\left( \begin{array}{c}
{\bf e}_{1}\\
{\bf e}_{2}  \\
{\bf e}_{3}
\end{array} \right)_{\xi_2} = C
\left( \begin{array}{c}
{\bf e}_{1}\\
{\bf e}_{2}  \\
{\bf e}_{3}
\end{array} \right)_{\xi_4} + D
\left( \begin{array}{c}
{\bf e}_{1}\\
{\bf e}_{2}  \\
{\bf e}_{3}
\end{array} \right).
\eqno(66)
$$

Hence we have the 3-dimensional M-LXX equation
$$
-bB_{\xi_{4}}+B_{\xi_{2}}-D_{\xi_{1}}+[B,D]=0
\eqno(67)
$$
which is the compatibility condition of the  equations (66) as $C=bI$.

\subsubsection{The 3-dimensional Self-Dual Yang-Mills equation}

To derive the SDYME in $d=3$ dimensions we consider the following
particular case of the M-LXVII equation [16]
$$
\left( \begin{array}{c}
{\bf e}_{1}\\
{\bf e}_{2}  \\
{\bf e}_{3}
\end{array} \right)_{\xi_1} =
(A_{1}-\lambda A_{3})
\left( \begin{array}{c}
{\bf e}_{1}\\
{\bf e}_{2}  \\
{\bf e}_{3}
\end{array} \right)
\eqno(68a)
$$
$$
\left( \begin{array}{c}
{\bf e}_{1}\\
{\bf e}_{2}  \\
{\bf e}_{3}
\end{array} \right)_{\xi_2} = \lambda
\left( \begin{array}{c}
{\bf e}_{1}\\
{\bf e}_{2}  \\
{\bf e}_{3}
\end{array} \right)_{\xi_4} + (A_{2}-\lambda A_{4})
\left( \begin{array}{c}
{\bf e}_{1}\\
{\bf e}_{2}  \\
{\bf e}_{3}
\end{array} \right).
\eqno(68b)
$$

The compatibility condition of these equations yields
 the 3-dimensional SDYME
$$
A_{2\xi_{1}}-A_{1\xi_{2}}+[A_{2},A_{1}]=0  \eqno(69a)
$$
$$
-A_{3\xi_{4}}+[A_{4},A_{3}]=0  \eqno(69b)
$$
$$
A_{1\xi_{4}}-A_{4\xi_{1}}+[A_{1},A_{4}]=
-A_{3\xi_{2}}+[A_{2},A_{3}].  \eqno(69c)
$$

\subsubsection{The Zakharov I equation}

Now let us we consider the M-LXVII equation in the form [16]
$$
\left( \begin{array}{c}
{\bf e}_{1}\\
{\bf e}_{2}  \\
{\bf e}_{3}
\end{array} \right)_{x} = (A_{1}-\lambda A_{3})
\left( \begin{array}{c}
{\bf e}_{1}\\
{\bf e}_{2}  \\
{\bf e}_{3}
\end{array} \right)
\eqno(70a)
$$
$$
\left( \begin{array}{c}
{\bf e}_{1}\\
{\bf e}_{2}  \\
{\bf e}_{3}
\end{array} \right)_{t} = \lambda
\left( \begin{array}{c}
{\bf e}_{1}\\
{\bf e}_{2}  \\
{\bf e}_{3}
\end{array} \right)_{y} + A_{2}
\left( \begin{array}{c}
{\bf e}_{1}\\
{\bf e}_{2}  \\
{\bf e}_{3}
\end{array} \right)
\eqno(70b)
$$
with
$$
A_{1}=
\left( \begin{array}{ccc}
0  &  i(q-p)  &  (q+p)\\
-i(q-p)  &  0  &  0  \\
-(q+p)  &  0  &  0
\end{array} \right)  \eqno(71a)
$$
$$
A_{2}=
\left( \begin{array}{ccc}
0  &  (q+p)_y  &  i(p-q)_y \\
-(q+p)_y  &  0  &  v  \\
-i(p-q)_y  &  -v  &  0
\end{array} \right)  \eqno(71b)
$$
$$
A_{3}=
\left( \begin{array}{ccc}
0  &    0  &  0 \\
0  &  0  &  1  \\
0  & -1  &  0
\end{array} \right).  \eqno(71c)
$$
The compatibility condition of these equations gives
$$
iq_{t}=q_{xy}+vq  \eqno(72a)
$$
$$
-ip_{t}=p_{xy}+vp  \eqno(72b)
$$
$$
v_{x}=2(pq)_{y}\eqno(72c)
$$
which for convenience we call the Zakharov I (Z-I) equation  [32].

\subsubsection{Integrable spin systems}

Now we consider the case when the curves/surfaces are given by
the following geometrical equation [16]
$$
\left( \begin{array}{c}
{\bf e}^{\prime}_{1}\\
{\bf e}^{\prime}_{2}  \\
{\bf e}^{\prime}_{3}
\end{array} \right)_{x} = -\lambda A_{3}
\left( \begin{array}{c}
{\bf e}^{\prime}_{1}\\
{\bf e}^{\prime}_{2}  \\
{\bf e}^{\prime}_{3}
\end{array} \right)
\eqno(73a)
$$
$$
\left( \begin{array}{c}
{\bf e}^{\prime}_{1}\\
{\bf e}^{\prime}_{2}  \\
{\bf e}^{\prime}_{3}
\end{array} \right)_{t} = \lambda
\left( \begin{array}{c}
{\bf e}^{\prime}_{1}\\
{\bf e}^{\prime}_{2}  \\
{\bf e}^{\prime}_{3}
\end{array} \right)_{y} -\lambda  A_{4}
\left( \begin{array}{c}
{\bf e}^{\prime}_{1}\\
{\bf e}^{\prime}_{2}  \\
{\bf e}^{\prime}_{3}
\end{array} \right)
\eqno(73b)
$$
with
$$
A_{3}=
\left( \begin{array}{ccc}
0  &  rS_1  & -irS_2 \\
-rS_1  &  0  &  S_3  \\
irS_2  &  -S_3  &  0
\end{array} \right), \quad A_{4}=
$$
$$
\left( \begin{array}{ccc}
0  & -ir[2iS_3 S_{2y}-2iS_2 S_{3y}+iuS_1] & -r[2S_3 S_{1y}-2S_1 S_{3y}-uS_2] \\
ir[2iS_3 S_{2y}-2iS_2 S_{3y}+iuS_1]  &  0  & -[ir^2(S^{+}S^{-}_{y}-S^{-}S^{+}_y)-uS_3]  \\
r[2S_3 S_{1y}-2S_1 S_{3y}-uS_2   & ir^2(S^{+}S^{-}_{y}-S^{-}S^{+}_{y})-uS_3   &  0
\end{array} \right).  \eqno(74)
$$
Here we have the additional condition
$$
S_{3}^{2}+r^{2}(S^{2}_{1}+S^{2}_{2})=1, \quad r^2=\pm 1. \eqno(75)
$$
The equation (73) with the condition (75) we call the M-LXVI equation.
From the compatibility condition of the equations (73) we obtain
the Myrzakulov I (M-I) equation
$$
iS_{t}= ([S,S_{y}]+2iuS)_{x}   \eqno(76a)
$$
$$
u_{x}=-\frac{1}{2i}tr(S[S_{x},S_{y}]) \eqno(76b)
$$
where
$$
S=
\left( \begin{array}{cc}
S_{3}  & rS^{-} \\
rS^{+}  &  -S_{3}
\end{array} \right), \quad S^{\pm}=S_{1}\pm iS_{2}. \eqno(77)
$$

\subsubsection{The (2+1)-dimensional mKdV equation}.

Let us now we consider the following version of the M-LXVII
equation
$$
\left( \begin{array}{c}
{\bf e}_{1}\\
{\bf e}_{2}  \\
{\bf e}_{3}
\end{array} \right)_{x} = (A_{1}-\lambda A_{3})
\left( \begin{array}{c}
{\bf e}_{1}\\
{\bf e}_{2}  \\
{\bf e}_{3}
\end{array} \right)
\eqno(78a)
$$
$$
\left( \begin{array}{c}
{\bf e}_{1}\\
{\bf e}_{2}  \\
{\bf e}_{3}
\end{array} \right)_{t} = \lambda^{2}
\left( \begin{array}{c}
{\bf e}_{1}\\
{\bf e}_{2}  \\
{\bf e}_{3}
\end{array} \right)_{y} + (D_{1}\lambda +D_{0})
\left( \begin{array}{c}
{\bf e}_{1}\\
{\bf e}_{2}  \\
{\bf e}_{3}
\end{array} \right)
\eqno(78b)
$$
where $A_{1}, A_{3}$ are given by (71) and $D_k $ are some matrices [16].
Then the complex functions $q,p$
satisfy the (2+1)-dimensional complex mKdV equation
$$
q_{t} +q_{xxy}-(qv_{1})_{x}-v_{2}q=0    \eqno(79a)
$$
$$
p_{t} +p_{xxy}-(pv_{1})_{x}-v_{2}p=0    \eqno(79b)
$$
$$
v_{1x}=2(pq)_{y}    \eqno(79c)
$$
$$
v_{2x}=2(pq_{xy}-p_{xy}q)   \eqno(79d)
$$
which is the 2+1 dimensional complex mKdV [30]. If $p=\beta q$ is real,
we get the following mKdV equation
$$
q_{t} +q_{xxy}-(qv_{1})_{x}=0    \eqno(80a)
$$
$$
v_{1x}=2\beta (q^2)_{y}.    \eqno(80b)
$$

\subsubsection{The (2+1)-dimensional derivative NLSE}.

Let  now we work with  the following form of the M-LXVII equation [16]
$$
\left( \begin{array}{c}
{\bf e}_{1}\\
{\bf e}_{2}  \\
{\bf e}_{3}
\end{array} \right)_{x} = (A_{3}\lambda^{2}+ A_{1}\lambda)
\left( \begin{array}{c}
{\bf e}_{1}\\
{\bf e}_{2}  \\
{\bf e}_{3}
\end{array} \right)
\eqno(81a)
$$
$$
\left( \begin{array}{c}
{\bf e}_{1}\\
{\bf e}_{2}  \\
{\bf e}_{3}
\end{array} \right)_{t} = \lambda^{2}
\left( \begin{array}{c}
{\bf e}_{1}\\
{\bf e}_{2}  \\
{\bf e}_{3}
\end{array} \right)_{y} + (D_{2}\lambda^{2}+D_{1}\lambda +D_{0})
\left( \begin{array}{c}
{\bf e}_{1}\\
{\bf e}_{2}  \\
{\bf e}_{3}
\end{array} \right)
\eqno(81b)
$$
where $A_{1}, A_{3}$ are given by (71).
Then the complex functions $q,p$
satisfy the Strachan equation [30]
$$
iq_{t}=q_{xy}-2ic(vq)_x  \eqno(82a)
$$
$$
-ip_{t}=p_{xy}+2ic(vq)_x  \eqno(82b)
$$
$$
v_{x}=2(pq)_{y}.  \eqno(82c)
$$

\subsubsection{The M-III$_{q}$ equation}.

At last  we consider the case
$$
\left( \begin{array}{c}
{\bf e}_{1}\\
{\bf e}_{2}  \\
{\bf e}_{3}
\end{array} \right)_{x} = (A_{3}(c\lambda^{2}+d\lambda)+ A_{1}(2c\lambda+d))
\left( \begin{array}{c}
{\bf e}_{1}\\
{\bf e}_{2}  \\
{\bf e}_{3}
\end{array} \right)
\eqno(83a)
$$
$$
\left( \begin{array}{c}
{\bf e}_{1}\\
{\bf e}_{2}  \\
{\bf e}_{3}
\end{array} \right)_{t} = 2(c\lambda^{2}+d\lambda)
\left( \begin{array}{c}
{\bf e}_{1}\\
{\bf e}_{2}  \\
{\bf e}_{3}
\end{array} \right)_{y} + (D_{2}\lambda^{2}+D_{1}\lambda +D_{0})
\left( \begin{array}{c}
{\bf e}_{1}\\
{\bf e}_{2}  \\
{\bf e}_{3}
\end{array} \right)
\eqno(83b)
$$
where $A_{1}, A_{3}$ are given by (71).
Then the complex functions $q,p$
satisfy the (2+1)-dimensional M-III$_{q}$ equation [16]
$$
iq_{t}=q_{xy}-2ic(vq)_{x}+d^{2}vq  \eqno(84a)
$$
$$
-ip_{t}=p_{xy}+2ic(vq)_{x}+d^{2}vp  \eqno(84b)
$$
$$
v_{x}=2(pq)_{y}.  \eqno(84c)
$$
The M-III$_{q}$ equation (84) admits two integrable reductions: the
Strachan equation (82) as $d=0$ and the Z-I equation (72) as $c=0$.

\section{Soliton geometry in $d=4$ dimensions}

There exist several equations of the soliton geometry in $d=4$ dimesions.
Some of them we present here.

\subsection{The M-LXVIII equation}
Consider the M-LXVIII equation [16]
$$
\left( \begin{array}{c}
{\bf e}_{1}\\
{\bf e}_{2}  \\
{\bf e}_{3}
\end{array} \right)_{\xi_1} = A
\left( \begin{array}{c}
{\bf e}_{1}\\
{\bf e}_{2}  \\
{\bf e}_{3}
\end{array} \right)_{\xi_3} + B
\left( \begin{array}{c}
{\bf e}_{1}\\
{\bf e}_{2}  \\
{\bf e}_{3}
\end{array} \right)
\eqno(85a)
$$
$$
\left( \begin{array}{c}
{\bf e}_{1}\\
{\bf e}_{2}  \\
{\bf e}_{3}
\end{array} \right)_{\xi_2} = C
\left( \begin{array}{c}
{\bf e}_{1}\\
{\bf e}_{2}  \\
{\bf e}_{3}
\end{array} \right)_{\xi_4} + D
\left( \begin{array}{c}
{\bf e}_{1}\\
{\bf e}_{2}  \\
{\bf e}_{3}
\end{array} \right)
\eqno(85b)
$$
where ${\bf e}^{2}_{j}=1, ({\bf e}_{i} {\bf e}_{j})=\delta_{ij}$ and
$A(\lambda), B(\lambda), C(\lambda), D(\lambda)$ are (3$\times$3)-matrices,
$\lambda$ is some parameter, $\xi_{i}$ are coordinates. This equation
describes some four dimensional curves and/or "surfaces" in 3-dimensional
space. It is one of main equations of the multidimensional soliton geometry
and admits several integrable reductions.

\subsection{The M-LXXI equation}

Consider the M-LXXI equation [16]
$$
\left( \begin{array}{c}
{\bf e}_{1}\\
{\bf e}_{2}  \\
{\bf e}_{3}
\end{array} \right)_{\xi_1} = A
\left( \begin{array}{c}
{\bf e}_{1}\\
{\bf e}_{2}  \\
{\bf e}_{3}
\end{array} \right), \quad
\left( \begin{array}{c}
{\bf e}_{1}\\
{\bf e}_{2}  \\
{\bf e}_{3}
\end{array} \right)_{\xi_{2}} =  B
\left( \begin{array}{c}
{\bf e}_{1}\\
{\bf e}_{2}  \\
{\bf e}_{3}
\end{array} \right)
\eqno(86a)
$$
$$
\left( \begin{array}{c}
{\bf e}_{1}\\
{\bf e}_{2}  \\
{\bf e}_{3}
\end{array} \right)_{\xi_{4}} = C
\left( \begin{array}{c}
{\bf e}_{1}\\
{\bf e}_{2}  \\
{\bf e}_{3}
\end{array} \right)_{\xi_{3}} + D
\left( \begin{array}{c}
{\bf e}_{1}\\
{\bf e}_{2}  \\
{\bf e}_{3}
\end{array} \right)
\eqno(86b)
$$

The compatibility condition of these equations gives some nonlinear evolution
equations (NEEs).

\subsection{The M-LXI equation}
Consider the 4-dimensional M-LXI equation [16]
$$
\left( \begin{array}{c}
{\bf e}_{1}\\
{\bf e}_{2}  \\
{\bf e}_{3}
\end{array} \right)_{\xi_1} = A
\left( \begin{array}{c}
{\bf e}_{1}\\
{\bf e}_{2}  \\
{\bf e}_{3}
\end{array} \right), \quad
\left( \begin{array}{c}
{\bf e}_{1}\\
{\bf e}_{2}  \\
{\bf e}_{3}
\end{array} \right)_{\xi_{2}} =  B
\left( \begin{array}{c}
{\bf e}_{1}\\
{\bf e}_{2}  \\
{\bf e}_{3}
\end{array} \right)
$$
$$
\left( \begin{array}{c}
{\bf e}_{1}\\
{\bf e}_{2}  \\
{\bf e}_{3}
\end{array} \right)_{\xi_{3}} = C
\left( \begin{array}{c}
{\bf e}_{1}\\
{\bf e}_{2}  \\
{\bf e}_{3}
\end{array} \right), \quad
\left( \begin{array}{c}
{\bf e}_{1}\\
{\bf e}_{2}  \\
{\bf e}_{3}
\end{array} \right)_{\xi_{4}} = D
\left( \begin{array}{c}
{\bf e}_{1}\\
{\bf e}_{2}  \\
{\bf e}_{3}
\end{array} \right).
\eqno(87)
$$
The compatibility condition of these equations gives the following
4-dimensional M-LXII equation
$$
A_{\xi_{2}}-B_{\xi_{1}}+[A,B]=0, \quad
A_{\xi_{3}}-C_{\xi_{1}}+[A,C]=0,\quad
A_{\xi_{4}}-D_{\xi_{1}}+[A,D]=0   \eqno(88a)
$$
$$
C_{\xi_{2}}-B_{\xi_{3}}+[C,B]=0, \quad
D_{\xi_{2}}-B_{\xi_{4}}+[D,B]=0, \quad
C_{\xi_{4}}-D_{\xi_{3}}+[C,D]=0.   \eqno(88b)
$$

This equation contents many interesting 4-dimnsional NEEs.

\subsection{The M-LXX equation}

From (85) we get the following M-LXX equation [16]
$$
AD_{\xi_{3}}-CB_{\xi_{4}}+B_{\xi_{2}}-D_{\xi_{1}}+[B,D]=0
\eqno(89a)
$$
$$
A_{\xi_{2}}-CA_{\xi_{4}}+[A,D]=0
\eqno(89b)
$$
$$
[A,C]=0
\eqno(89c)
$$
$$
C_{\xi_{1}}-AC_{\xi_{3}}+[C,B] =0.
\eqno(89d)
$$

If we choose
$$
A=aI, \quad C=bI, \quad a,b =consts  \eqno(90)
$$
then the M-LXX equation (89) takes the form
$$
aD_{\xi_{3}}-bB_{\xi_{4}}+B_{\xi_{2}}-D_{\xi_{1}}+[B,D]=0.
\eqno(91)
$$

\subsection{The SDYME }

Now we assume that
$$
B=A_{1}-\lambda A_{3}, \quad D=A_{2}-\lambda A_{4}, \quad a=b=\lambda.
\eqno(92)
$$
So that the M-LXVIII equation takes the form [16]
$$
\left( \begin{array}{c}
{\bf e}_{1}\\
{\bf e}_{2}  \\
{\bf e}_{3}
\end{array} \right)_{\xi_1} = \lambda
\left( \begin{array}{c}
{\bf e}_{1}\\
{\bf e}_{2}  \\
{\bf e}_{3}
\end{array} \right)_{\xi_3} + (A_{1}-\lambda A_{3})
\left( \begin{array}{c}
{\bf e}_{1}\\
{\bf e}_{2}  \\
{\bf e}_{3}
\end{array} \right)
\eqno(93a)
$$
$$
\left( \begin{array}{c}
{\bf e}_{1}\\
{\bf e}_{2}  \\
{\bf e}_{3}
\end{array} \right)_{\xi_2} = \lambda
\left( \begin{array}{c}
{\bf e}_{1}\\
{\bf e}_{2}  \\
{\bf e}_{3}
\end{array} \right)_{\xi_4} + (A_{2}-\lambda A_{4})
\left( \begin{array}{c}
{\bf e}_{1}\\
{\bf e}_{2}  \\
{\bf e}_{3}
\end{array} \right).
\eqno(93b)
$$

From (91) we obtain the SDYME
$$
A_{2\xi_{1}}-A_{1\xi_{2}}+[A_{2},A_{1}]=0  \eqno(94a)
$$
$$
A_{4\xi_{3}}-A_{3\xi_{4}}+[A_{4},A_{3}]=0  \eqno(94b)
$$
$$
A_{1\xi_{4}}-A_{4\xi_{1}}+[A_{1},A_{4}]=
A_{2\xi_{3}}-A_{3\xi_{2}}+[A_{2},A_{3}]  \eqno(94c)
$$
or
$$
F_{\xi_{1}\xi_{2}}=0, \quad F_{\xi_{3}\xi_{4}}=0, \quad
F_{\xi_{4}\xi_{1}}-F_{\xi_{3}\xi_{2}}=0.  \eqno(94d)
$$
Here
$$
F_{\xi_{i}\xi_{k}} = A_{k\xi_{i}}-A_{i\xi_{k}}+[A_{k},A_{i}].
$$
The SDYME (94) on a connection $A$ are the self-duality conditions
on the curvature under the Hodge star operation
$$
F= \ast F   \eqno(95a)
$$
or in index notation
$$
F_{\mu\nu}= \frac{1}{2}\epsilon_{\mu\nu\rho\delta}F^{\rho\delta}   \eqno(95b)
$$
where $\ast$ is Hodge operator, $\epsilon_{\mu\nu\rho\delta}$ stands
for the completely antisymmetric tensor in four dimensions with the convention:
$\epsilon_{1234}=1$. The SDYME is integrable by the Inverse Scattering
Transform  method (see, e.g. [28,29]). The Lax representation (LR) of the SDYME has the form [28, 29]
$$
\Phi_{\xi_{1}}-\lambda\Phi_{\xi_{3}}=(A_{1}-\lambda A_{3})\Phi  \eqno(96a)
$$
$$
\Phi_{\xi_{2}}-\lambda\Phi_{\xi_{4}}=(A_{2}-\lambda A_{4})\Phi.  \eqno(96b)
$$
Hence follows that for the SDYME the spectral parameter $\lambda$ satisfies
the equations
$$
\lambda_{\xi_{1}}=\lambda\lambda_{\xi_{3}}, \quad
\lambda_{\xi_{2}}=\lambda\lambda_{\xi_{4}}. \eqno(97)
$$
These equations have the following solutions
$$
\lambda=\frac{n_{1}\xi_{3}+n_{3}}{n_{4}-n_{1}\xi_{1}}, \quad
\lambda=\frac{m_{1}\xi_{4}+m_{3}}{m_{4}-m_{1}\xi_{2}}.  \eqno(98)
$$
So that the general solution of the set (97) has the form
$$
\lambda=\frac{n_{1}\xi_{3}+n_{3}+m_{1}\xi_{4}}{n_{4}-n_{1}\xi_{1}-
m_{1}\xi_{2}}  \eqno(99)
$$
where $m_{i}, n_{i} =constants$. The corresponding solution of the SDYME (94)
is called the breaking (overlapping) solutions.

\section{Conclusion}
In this note, we have formulated the some classes of the
motion of curves/surfaces
in $d=3, 4$-dimensional space using the differential geometry. It is shown
that some of these curves/surfaces are integrable in the sense that
they are connected with the well known integrable equations
in multidimensions.  Examples include practically all
known multidimensional integrable (soliton) equations:
the DS, Zakharovs, NLS-types, (2+1)-KdV,  mKdV , Ishimori, M-IX,
M-I   equations and the SDYME. In particular,  we have shown that one of
(2+1)-dimensional extensions of the GMCE (10), namely, {\bf the M-LXII (or
mM-LXII) equation (19)=(59) is integrable} as the particular case of the
integrable BE (20). It means that {\bf the M-LXII equation is exact
reduction of the famous SDYME}. In turn, it indicate that  {\bf almost all
known soliton eqiations in 2+1 dimensions are exact reductions of the SDYME}
since these equations obtained from the mM-LXII (or M-LXII) equation as some
particular cases (see for instance refs [?-?]).

Although main elements of our approach have been  established,
but there remain  many problems to be studied. The study of  some
of these  problems will be the subjects of our   future  works. Here only
we note that there exists the  other approach to study of the
multidimensional
soliton geometry developing mainly by Konopelchenko and coworkers (see, e.g.
refs. [9, 11-14] and references therein). The main tool of this
approach is a generalized Weierstrass representation for a conformal
immersion of surfaces into $R^{3}$ or $R^{4}$, and also a linear
problem related with this representation. A consideration of the
linear  problem along with the Weierstrass representation
allows to express integrable deformations of surfaces via such hierarchies
soliton equations as a Veselov-Novikov hierarchy,   DS hierarchy and so on.
Thus, in this context,  this and our approaches, developing in paralleel,
have the
common  purpose, namely, the construction surfaces (and curves) inducing
by multidimensional soliton equations and  complement to each other.

\section{Acknowledgments}

RM would like to thanks to Prof. M.Lakshmanan
for stimuliting discussions, hospitality during visits and for the
financial support.

\end{document}